\documentclass[a4paper,12pt]{article}
\usepackage[utf8]{inputenc}
\usepackage[english]{babel}
\usepackage{amssymb}
\usepackage{amsmath}
\usepackage{epsfig}
\usepackage{bm}
\usepackage{graphics}\textwidth=17cm
\def\be{\begin{equation}}
\def\ee{\end{equation}}
\def\ord{{\rm ord}}
\def\deg{{\rm deg}}

\newtheorem{thrm}{\bf Theorem}

\graphicspath{{../img/}}

\begin{document}

\title{Stiefel Filters}
\author{Andrei Bogatyr\"ev
\thanks{Supported by INM RAS branch of MCFAM, agreement 075-15-2022-286.}
}
\date{}
\maketitle

\abstract{The best uniform rational approximation of the \emph{Sign} function on two
intervals separated by zero was explicitly found by E.I. Zolotar\"ev in 1877. 
The natural extension of this problem to three bands was solved by E.Stiefel in 1961.
We indicate  the solutions overlooked by the prominent geometer and study their properties.}

\noindent
{\bf Keywords:} {Uniform rational approximation, optimization of electrical filters, Ansatz method, equiripple property, alternation, Stiefel class}\\
{\bf MSC2010:} 41A20, 41A50, 49K35, 94Cxx\\

\section*{Introduction}
An electrical filter is a device which processes the harmonic components of the input signal in accordance to the transfer function. A signal component  is either supressed when its frequency lies in the given stopbands or it is translated without significant change in its magnitude when the frequency belongs to the given passbands. The optimal synthesis of such devices is based on a uniform rational approximation problem which is a generalization of the renowned  third (or fourth) Zolotar\"ev problem \cite{Zol}. Roughly, the problem consists in the best uniform approximation of the two-valued (indicator) function defined on the prescribed  passbands and stopbands -- the collection $E$ of disjoint segments on real frequency axis -- by a rational function of bounded degree $n$. 

The simplest problem of this type -- for two bands --  was completely solved by Chebyshev's pupil E.I.Zolotar\"ev yet in 1877 \cite{Zol}. Ed Stiefel in 1961 made the next natural step and gave a  solution to the problem with three bands \cite{Sti}. An algebro-geometric Ansatz solving the general filter approximation problem \cite{B10} was proposed  in 2010. The latter solution has the following structure: it is the elliptic sine of an abelian integral and generalizes the representation for Zolotar\"ev fractions. It contains the unknown parameters related to the Riemann surface that bears the abelian integral, both continuous and discrete. The latter have to be evaluated given the input data of the problem. The number of those parameters is usually greatly less than the number of parameters in the rational function we optimize, 
which gives us the considerable dimensional reduction of the problem.

In this paper we analyse the three band case from the viewpoint of the Ansatz method which has the profound links to the theory of Belyi functions and Grothendieck's Dessins d'Enfants \cite{Gr, LZ} and their generalizations. The reason for that is the following observation.
A degree $n$ solution of the uniform approximation problem has $2n+2$ so called alternation points, each of those in the interiour of the work bands is necessarily a critical point of the solution with the value in the four element set defined by the error of approximation.  Their number almost coinsides with the total number $2n-2$ of critical points of a degree $n$ rational function. For instance, Zolotarev fraction has exactly four distinct critical values and this property distinguishes it among other rational functions \cite{B17}. Compare it to characteristic property of Chebyshev polynomials which  have just two separate finite critical values. The graphs which are the counterparts of Dessin d'Enfants for this theory completely describe the discrete part of the Ansatz and may be explicitly listed manualy when the number of bands is not too big (or by a computer otherwise). The job was done for three bands and showed that Stiefel's formula has to be modified for certain sets of bands. In particular, the genus $g$ of the  associated algebraic curve (which was assumed to be two by Stiefel) varies from $g=1$ to $g=3$ for the general set of three bands. We may suppose that the case $g=2$ is dominant in some sense, however two other cases cannot be considered as exceptional: each of those correspond to an open set in the space of three band problems.  

All statements announced without a proof in this paper will be considered in greater detail in the extended version of this work.

\section{Optimization problem for a multiband filter}
Suppose a finite collection $E$ of disjoint closed segments of extended real line $\hat{\mathbb{R}}=\mathbb{RP}^1$ is given. The set has a meaning of work frequency bands of a filter, its compliment on the frequency axis is called the set  $T$ of  transient bands. Work bands are of   two types: $E=E^+\cup E^-$,  the passbands $E^+$ and the stopbands $E^-$. Both subsets $E^\pm$ are non empty. Optimization problem for electrical filter has several equivalent settings \cite{Cauer, AS, Akh, B21}. Here we consider just one of them -- the natural multiband extension of the fourth Zolotarev problem \cite{Ach2}.

Define the indicator function $S_E(x)=\pm1$ when $x\in E^\pm$. Find the best uniform  approximation $R(x)$  of the function $S_E(x)$  
\be
\label{Appr}
 ||R-S_E||_{C(E)} := \max\limits_{x \in E}|R(x) - S_E(x)| \to \min =: \mu.
\ee
among all real rational functions $R(x)$ of bounded topological degree $\deg R\le n$ (being the maximum of the  degrees of the numerator and the denominator for the non-cancellable fraction). 

In this paper we consider three  band filter  problems only 
and assume w.l.o.g. that there is one stopband and two passbands.
The following decomposition of the extended real line into cyclically ordered nonintersecting work ($E$) and transition ($T$) bands 
\be
E^-<T_1<E_1^+<T_{12}<E_2^+<T_2<E^-.
\label{AllBands}
\ee
gives standard labels to the latter.

\subsection{Some remarks}
 Already the first investigations \cite{Zol,Akh} of the approximation problem \eqref{Appr} showed that even in the classical Zolotar\"ev case with one pass- and one stop- band,  the goal function \eqref{Appr} has local minima. This phenomenon was observed in the Stiefel's paper \cite{Sti} and fully explained in the  dissertation of his pupil R.-A.R.Amer \cite{AS} who decomposed the space of rational functions of bounded deviation into classes. In the closure of each nonempty class there is a unique solution of the above optimization problem \cite{B21}. 
  
  In the current problem setting we define the topological class $\sigma$ of a competing rational function $R(x)$ with deviation $\mu<1$ by fixing the parity of the number of its zeros (alternatively: poles) in each transient band. For the three band problem $\sigma(R)=(\sigma_1,\sigma_{12}, \sigma_2)\in\mathbb{Z}_2^3$ where $\sigma_*$ is related to the transient band $T_*$, $*=1,12,2$. There is an obvious constraint: $\sigma_1+\sigma_{12}+\sigma_2=\deg R ~mod~2$.

   All local optimal functions are characterized by the \emph{alternation} (\emph{equiripple} in terms of electrical engineers) property \cite{Ach2, B21}. In our problem setting the approximation error $\delta(x):=R(x)-S_E(x)$ of degree $n$ minimizer has $2n+2$ \emph{alternation} points $a_s\in E$ where  $\delta(a_s)=\pm||\delta||_{C(E)}$ with the consecutive change of sign. For the solution with positive defect $n-\deg ~R$ the criterion of optimality is also known but it is more sophisticated \cite{B21}. Roughly, the equiripple property says that the graph of (the error function of) a solution on the workbands $E$ looks like the sequence of sufficiently many ripples of constant amplitude. With this property taken into account, the solution of the problem becomes in a sense very simple: you merely produce  a function with the requested  oscillatory behaviour. This is our workplan for the next three sections of the paper.

\section{Two band case}\label{sect:Zol}
We start with the brief description of Zolotarev fraction construction. The latter may be turned into the Stiefel's solution by excitation of some internal degree of freedom. Zolotarev's original work \cite{Zol} cannot be called transparent: it is based on the modular theory of elliptic functions which additionally is treated as a purely real theory (following the ideology of  P.L.Chebyshev). N.I.Achiezer has significantly simplified solution of Zolotarev re-writing it in the adequate language of complex variables \cite{Akh1, Ach2}. Still it remains rather difficult for the understanding. Following \cite{B12} we give a clear geometric picture which easily explains all features of Zolotarev fraction.

Let $\Pi(\tau)$, $\tau\in i(0,\infty) $, be a rectangle of dimensions $2\times|\tau|$:
\be
\label{Rectangle}
\Pi(\tau):=\{u\in \mathbb{C}:\quad -1< Re~u< 1,~~ 0< Im~ u< |\tau| \}.
\ee
The rectangle may be conformally mapped 1-1  to the upper half plane fixing three points on the boundary:
\be
x(u|\tau):\quad  \Pi(\tau),~-1,~0,~1 \to \mathbb{H},~-1,~0,~1. 
\ee
This mapping is computationally efficient, it is given essentially by the elliptic sine of modulus $\tau$ and may be expressed via standard elliptic thetas \cite{Akh1} implicitly depending on the modulus $\tau$:
\be
x(u|\tau)= sn(K(\tau)u|\tau)=\frac{\theta_3}{\theta_2}~\frac{\theta_1(u/2)}{\theta_0(u/2)},
\ee
where $sn(\cdot|\tau)$  is the elliptic sine of modulus $\tau$, $K(\tau)=\frac\pi2\theta_3^2$ is the complete elliptic integral; and theta without argument is the theta constant i.e. the value for the argument $u=0$.

\begin{figure}
\includegraphics[width=\textwidth, trim=0 24.5cm 0 0, clip]{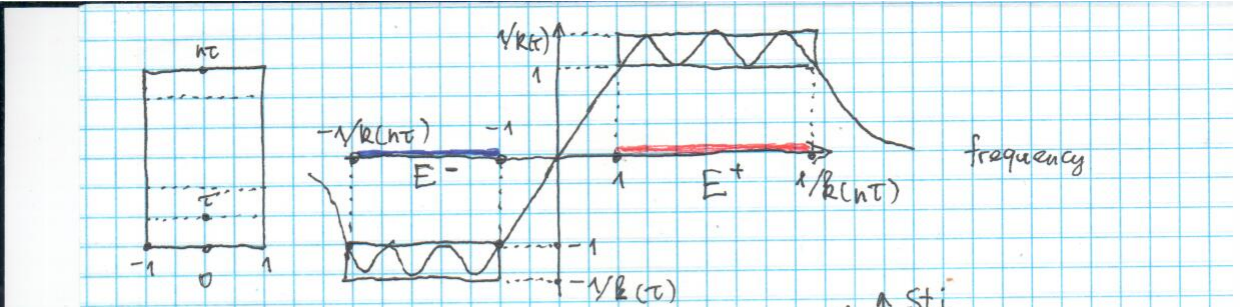}
\caption{Left:~~ Big rectangle tiled by the smaller ones, $n=6$\\
Right:~~ The qualitative graph of Zolotarev fraction.}
\label{ZolGraph}
\end{figure}

The rectangle $\Pi(n\tau)$ is naturally tiled by $n>0$ copies of smaller rectangle $\Pi(\tau)$ as shown in the Fig. \ref{ZolGraph}. The conformal developing maps of the bigger and  the smaller rectangles are related by a \emph{rational} degree $n$ map $Z_n$:
\be
Z_n(x(u|n\tau))=x(u|\tau).
\ee
This statement easily follows from the Schwarz reflection principle for the conformal mapping $x(\cdot|\tau)\circ x^{-1}(\cdot|n\tau)$ sending upper half plane to the Riemann sphere.
Indeed, this map takes real values at the boundary of the half plane and hence admits an analytic continuation to the whole Riemann sphere where it has exactly $n$ simple poles as the only singularities. The poles correspond to points $u\in (2\mathbb{Z}+1)\tau$.

It is not difficult to draw the qualitative graph of this function on the real axis:
$Z_n$ is monotonic outside two segments $E^\pm:=\pm[1,1/k(n\tau)]$ $=x(u|n\tau)$, $u\in \pm1+[0,n]\tau$, and oscillates exactly $n$ times on each of them between the values $\pm\{1,1/k(\tau)\}$ (see the right panel of Fig. \ref{ZolGraph}). 
Here $1/k(\tau):=x(1+\tau|\tau)=(\theta_3/\theta_2)^2$ is the value of an elliptic modular function. Therefore, a rescaled function $\frac{2k}{k+1}Z_n(x) $ has exactly $2n+2$ alternation points with the deviation value $\mu=(1-k)/(1+k)$  on the set of two bands $E^\pm$ which correspond to all corners of $n$ smaller rectangles. 

Zolotarev fractions possess lots of curious properties and often appear as the solutions 
of extremal problems of geometric function theory \cite{Dubini21, Dubini22}. In particular, Chebyshev polynomials are the special limit case of the suitably renormalized Zolotarev fractions when $\tau\to 0$. For instance, any Zolotarev fraction has exactly four different critical values: $\{\pm1,\pm1/k(\tau)\}$ which is essentially their characteristic property \cite{B17}; the composition of two fractions with suitably related modules is again a Zolotarev fraction \cite{B12} etc.

\section{Stiefel function as a deformation of Zolotarev fraction} 
\label{Sect:Sti}
Consider the large rectangle $\Pi(n\tau)$ with a slit at a quantized height, namely  along the boundary of one of the tiles:
\be
\label{Pi0}
\Pi(\tau; n,m; h):= \Pi(n\tau)\setminus\{[h,1]+m\tau\},
\qquad m=1,2,\dots,n-1;
\qquad  -1<h<1.
\ee
 We can conformally map the slit rectangle to the upper half-plane by a function\footnote{We share the concept of polymorphism: functions are destinguished by their names AND the 
 number/type of arguments.}
  $x(u|\tau, n, m, h)$. It will be convenient not to normalize this map, so the function is defined up to left action of $PSL_2(\mathbb{R})$ represented by real linear fractional maps conserving the orientation of the real line. The argument with the reflection principle  for the Zolotarev fraction remains valid, so the following Stiefel filtering function mapping the upper half plane to the slit rectangle and then to the Riemann sphere  
 \be
 St(\cdot|\tau;n,m;h):= x(\cdot|\tau)\circ x^{-1}(\cdot|\tau;n,m;h). 
 \ee
 is a real rational one of degree $n$.
 
 \begin{figure}
 \includegraphics[width=\textwidth, trim=0 18.5cm 0 5cm, clip]{ZS012.pdf}
  \caption{Left: ~Tiled rectangle with a slit, $n=7, m=4$.\\
 Right: Qualitative graph of the function $St(x|\dots)$. Bounding boxes of oscillations 
 project to the (maximal) work bands of the filter.}
 \label{StiGraph}
 \end{figure}

  One can draw the qualitative graph of this function -- see the left panel of Fig \ref{StiGraph}. We observe that after the multiplication of $St(x|\dots)$ by the factor of $2k/(k+1)$, $k=k(\tau)$, the function acquires $2n+2$ alternation points   on a set $E$ defined as follows. The stopband $E^-$ is the image of the left side of the slit rectangle \eqref{Pi0} under the map $x(u|\tau;n,m;h)$,  $u\in -1+[0,n\tau]$.  The choice of two passbands $E^+$ is twofold: they are the image of the right side of the slit rectangle with removed segment of the length not greater that $\tau$ and either beginning or ending at $1+m\tau$ -- the end point of the slit:
 
 \be
 E^+(v):=x(u|\tau;n,m;h), u\in\left\{
\begin{array}{l}
 1+[0,m]\tau ~or~ 1+[v,n]\tau,\qquad m\le v\le m+1;\\
 1+[0,v]\tau ~or~ 1+[m,n]\tau,\qquad m-1\le v\le m.
\end{array}
\right.
\label{StiE+}
 \ee
 
 The topological class $\sigma$ of this solution is $(1,0,n+1)~mod~2.$
 
 The mapping of the upper half plane to the slit rectangle \eqref{Pi0} has the Christoffel-Schwarz integral representation, which in our setting becomes a genus 2 hyperelliptic integral. So we  obtain the (universal, see below) representation of a Stiefel solution as an elliptic sine of a  genus two holomorphic abelian integral.
 
Stiefel has obtained the formula for this optimal function essentially in the same way as Chebyshev: writing Pell type functional equation which reduces to an ODE admitting separation of variables.  His presentation of the solution has the apperance of the equality of the elliptic and ultraelliptic (i.e. genus 2) integrals with variable upper limits being the argument and the value of the function.   
 
Thus constucted function $St(x|\dots)$ depends on two continuous parameters $h,\tau$ but solves the three dimensional set of optimization problems for $E=E(h,v,\tau)$. It is a known effect in uniform approximation that a single function may serve as a solution to the variety of problems. This happens because the ball in the C-norm has corners, see the discussion in the Introduction to \cite {Bbook}. 
 
The entire set of three band filter problems is also three-dimensional: 
the endpoints of 3-band set $E$ depend on 6 real parameters, but two sets $E$ which differ by an element of three-dimensional group $PSL_2(\mathbb{R})$ give us essentially the same problem. Below we show that Stiefel' solution does not cover all possible three band problems. 

\section{Other similar constructions}
It will be shown below that for a three band filter problem the genus of the algebraic curve participating in the Chebyshev Ansatz varies from 1 to 3. 
Preliminary we describe the constructions of appropriate solutions. Genus 2 case had been just described.

\subsection{Genus one}
Take a Zolotarev fraction and remove a segment from its passband $E^+$ that does not contain an alternation point inside it. We get two passbands while the number of alternation points did not change, hence Zolotarev fraction is a solution
for some (full dimensional) set of three band problems.
\be
\begin{array}{ll}
E^-&=x(u|n\tau),\qquad  u\in-1+[0,n]\tau;\\
E^+(v_1,v_2)&=x(u|n\tau),\qquad  u\in 1+[0,v_1]\tau ~or~ 1+[v_2,n]\tau,
\end{array}
\ee
where $m\le v_1< v_2\le m+1$ for some $m=0,1,2\dots,n-1$. We also exclude cases $v_1=0(=m)$ and $v_2=n(=m+1)$ to avoid the collapse of a passband to a point. This solution  has topological class $\sigma=$ $(1,0,n+1)~mod~2$.

\subsection{Genus three}
Consider four rectangular polygons in the plane, three of which are self-overlapping (two- sheeted) including one with interior branching. For good, these three live on a suitable Riemann surfaces, not plane. 

\subsubsection{Univalent case}
\label{Sect:Unival}
Let us excite yet another internal degree of freedom in Zolotarev's construction.
Consider a rectangle with two horizontal slits separated by a 'quant' $\tau$ in vertical direction:
\be
\label{Pi1}
\Pi_1(\tau;~n,m;~h_1,h_2):=\Pi(n\tau)\setminus\left(\{m\tau+[h_1,1]\}\cup
\{(m+1)\tau+[h_2,1]\}\right)
\ee
here the parameters of slits are chosen as follows: $h_1,h_2\in (-1,1)$, $m=1,2,\dots,n-2$. This polygon is conformally mapped to the upper half plane by the function $x_1(u|\tau; ~n, m;~h_1,h_2)$ defined up to the postcomposition with an element of $PSL_2(\mathbb{R})$ realised as a linear-fractional function.
The reasoning with the reflection principle from Section \ref{sect:Zol}
holds true, so the composition 
\be
St_1(\cdot|\tau; n, m;h_1,h_2):= x(u|\tau)\circ x_1^{-1}(u|\tau; n, m;h_1,h_2)
\label{Sti1}
\ee 
is a degree $n$ rational function whith qualitative graph  shown in Fig. \ref{Sti1Graph}. The function \eqref{Sti1} oscillates $n$ times on the stopband $E^-$ being the image of the left  wall of the twice slit rectrangle \eqref{Pi1} under the  mapping $x_1(u|\dots)$. Also it oscillates  exactly $m$ (resp. $n-m-1$) times on the first (resp. second) passband $E_1^+$
(resp. $E^+_2$) which is the image of the right wall below (resp. above) all the slits: 
\be
\begin{array}{l}
E^-=x_1(u|\tau; n, m;h_1,h_2),\qquad  u\in-1+[0,n]\tau;\\
E_1^+=x_1(u|\tau; n, m;h_1,h_2),\qquad  u\in 1+[0,m]\tau;\\
E_2^+=x_1(u|\tau; n, m;h_1,h_2),\qquad  u\in 1+[m+1,n]\tau,
\end{array}
\ee
Note that essentially yet another passband appeares  in the transient band $T_{12}$ with just one oscillation on it. There are two 'irregular' critical points (that is away  from the alternation set) of the function $St_1$, one in each gap between the 'fake' passband $x_1(1+[m,m+1]\tau|\dots)$ and $E_1^\pm$: they correspond to the tips of the slits.

\begin{figure}[h]
 
 \includegraphics[width=\textwidth, trim=0 12.6cm 0 12cm, clip]{ZS012.pdf}
 \caption{Left:Tiled rectangle with two slits, $n=8, m=4$.\\
 Right: Qualitative graph of the function $St_1(x|\dots)$. Bounding boxes project to the work bands of the filter.}
 \label{Sti1Graph}
 \end{figure}

Again one can count, that this function (after suitable rescaling) has the required number of alternation points on the set $E=E^-\cup E_1^+\cup E_2^+$ of bands and is therefore a solution.
Its topological class $\sigma$ is $(1,0,n+1)~mod~2$.
This construction has three free parameters and we can assume that it gives a 
solution to a full dimensional set of three-band problems.

The inverse to the conformal map $x_1(u|\dots)$ has an integral (SC) representation 
which is a holomorphic integral on a curve of genus three  with the maximal number (i.e. four) of real ovals. Six branchpoints of the curve are fixed and coinside with the 
endpoints of the working bands $E$. Two other branchpoints are free and have to be found 
along with the modulus $\tau$ in the process of the solution. They are the parameters of the Chebyshev Ansatz and make up a solution to the set of three equations.  

%
%

\subsubsection{Octagon with interior branching}
\label{Sect:Octagon}
\begin{figure}[h]
 
 \includegraphics[width=\textwidth, trim=0 6.7cm 0 16.8cm, clip]{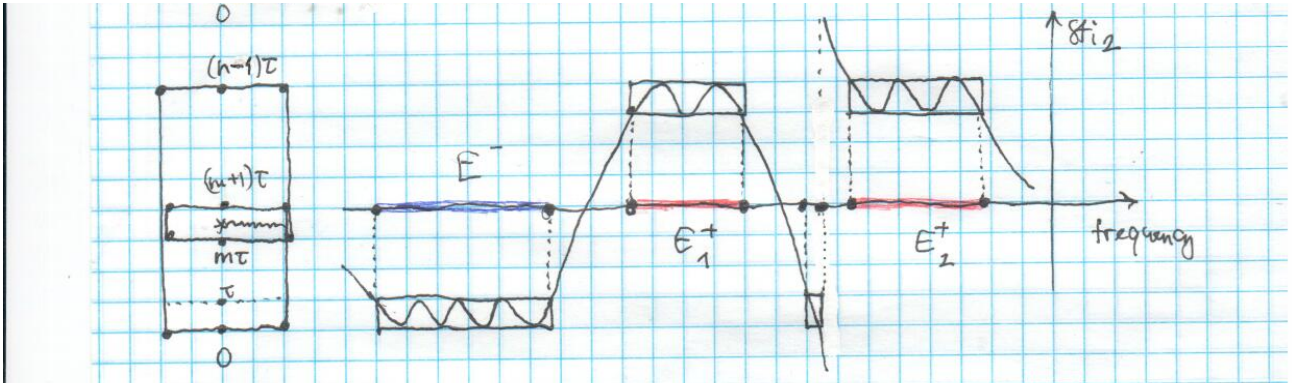}
 \caption{Left: Tiled rectangular octagon with interior branching $\ast$\\
 Right: Qualitative graph of the function $St_2(x|\dots)$. Bounding boxes of oscillations project to the work bands of the filter.}
 \label{Sti2Graph}
 \end{figure}

Let us consider the rectangular octagon with an inner branchpoint. 
\be
\label{Pi2}
\Pi_2(\tau;n,m;c):= \Pi((n-1)\tau)~~~ \leftarrow (c) \rightarrow ~~~\Pi(\tau)+m\tau,
\qquad m=1,2,\dots, n-2,
\ee
here $c\in\Pi(\tau)+m\tau,$ is the branchpoint and $\leftarrow(c)\rightarrow$ is the result of the following surgery operation sewing smaller rectangle to the larger one, see left panel of Fig. \ref{Sti2Graph}. Consider a slit in a smaller rectangle joining the branchpoint $c$ to a point at a common boundary
$1+\tau[m,m+1]$ of the rectangles. Draw the same slit in the larger rectangle and identify the right bank of one slit to the left bank of the other and vice versa. The result of the surgery is a topological disc with the conformal structure independent of the choice of the slit. This disc is  conformally mapped to the upper half plane by a function $x_2(u|\tau;n,m;c)$ (again, we do not specify the normalization of the map, so $PLS_2(\mathbb{R})$ acts on the left).  Our usual reasoning suggests that the function 
\be
\label{Sti2}
St_2(\cdot|\tau;n,m;c):= x(\cdot|\tau)\circ x_2^{-1}(\cdot|\tau;n,m;c). 
\ee
is a real rational function of degree $n$. Its qualitative graph is shown on the right panel of the figure \ref{Sti2Graph}. Note that $St_2$ has $n-1$ full oscillations 
at the stopband $E^-$  being the image of the longer left  wall of the  rectrangular octagon under the  mapping $x_2(u|\dots)$. Also it oscillates exactly $m+1$ and $n-m-1$ times on two  passbands $E_1^+$ and $E^+_2$ which are the images of the two right walls of the octagon: 
\be
\begin{array}{l}
E^-=x_2(u|\tau; n, m; c),\qquad  u\in-1+[0,n-1]\tau;\\
E_1^+=x_2(u|\tau; n, m;c),\qquad  u\in 1+[0,m+1]\tau;\\
E_2^+=x_2(u|\tau; n, m;c),\qquad  u\in 1+[m,n-1]\tau.
\end{array}
\ee

Totally, after proper rescaling  the function \eqref{Sti2} has the required number of alternation points on the aforementioned set $E$ of workbands and  gives us a solution of three-parametric set of three band filter problems. This solution  has topological class $\sigma=$ $(1,1,n)~mod~2$.

The intermediate mapping $x_2(u|\dots)$ of presentation \eqref{Sti2} is the inverse of the SC integral which in this particular case is the holomorphic abelian integral on a hyperelliptic genus 3 curve. Six branchpoints of this curve are the endpoints of the workbands $E$. Two more branchpoints are located in the gap between the passbands $E_1^+$, $E_2^+$ and bound the 'fake' stopband with oscillation number $1$. Along with the elliptic module $\tau$ which determines the deviation, those two free branchpoints should be determined from the input data by solution of three equations of the Ansatz method. 

%

 \subsubsection{Overlapping decagons}
 \label{Sect:Decagon}
\begin{figure}
\includegraphics[width=\textwidth, trim=0 14.5cm 0 .5cm, clip]{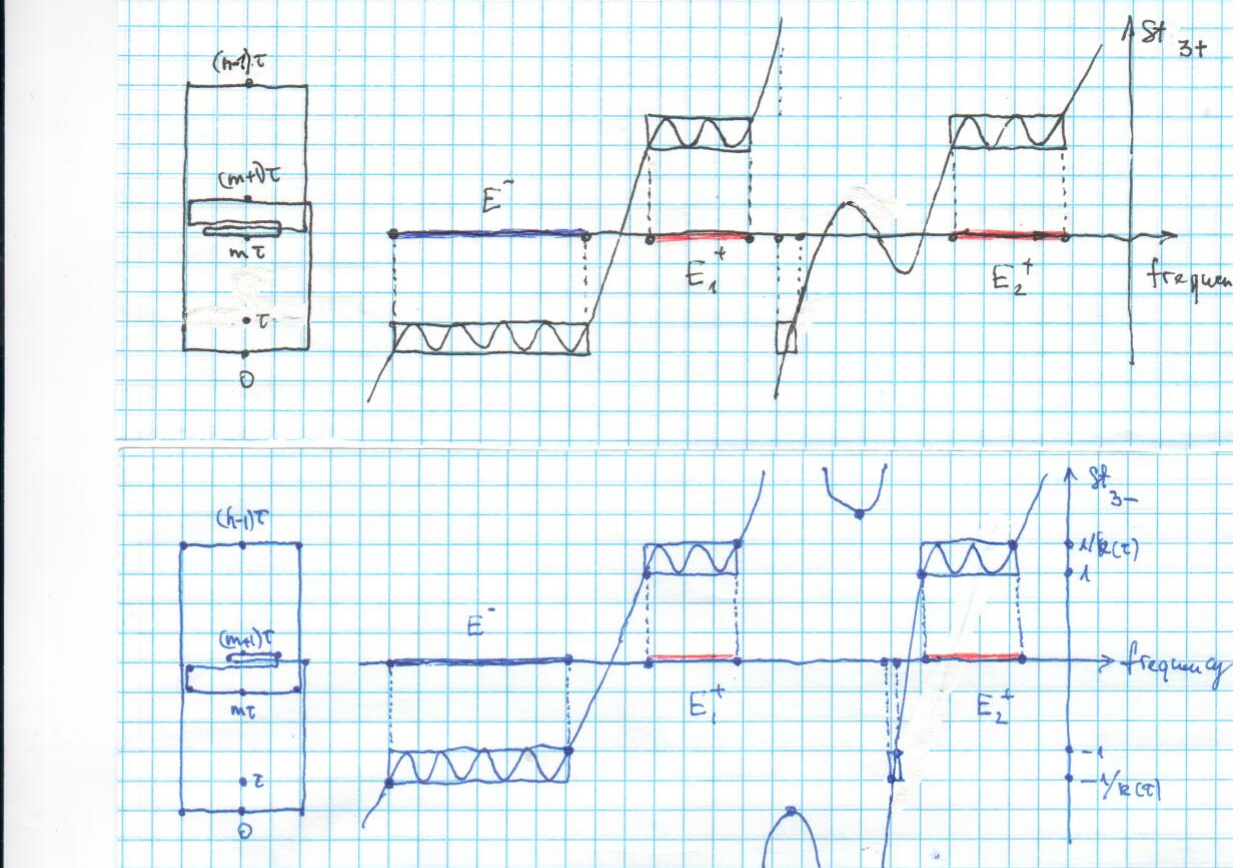} 
\caption{Left: Two rectangular decagons $\Pi_3^+$ (top), and $\Pi_3^-$ (bottom) for $n=10$ and $m=5$.\\
Right: Qualitative graphs of the functions $St_{3+}(x|\dots)$ (top) and $St_{3-}(x|\dots)$ (bottom). Bounding boxes project to the work bands of the filter}
 \label{Sti3Graph}
\end{figure}

 We consider two self- overlapping rectangular decagons which are topological cells. The tailoring kit is almost the same in both cases:
 \be
 \begin{array}{l}
 \Pi_3^+(\tau,n,m,h_1,h_2):= \Pi((n-1)\tau)\setminus\{[h_1,1]+m\tau\}~~~~
\leftarrow(h_2)\rightarrow ~~~~\Pi(\tau)+m\tau,\\
\Pi_3^-(\tau,n,m,h_1,h_2):= \Pi((n-1)\tau)\setminus\{[h_1,1]+(m+1)\tau\}~
\leftarrow(h_2)\rightarrow ~\Pi(\tau)+m\tau,
\end{array}
\ee
Here integer parameter $m$ ranges in $1,2,\dots,n-2$ for $\Pi_3^+$
and in $0,1,\dots,n-3$ for $\Pi_3^-$;  $-1<h_1\le h_2<1$. The glueing procedure $\leftarrow(h_2)\rightarrow$ also depends on the superscript $\pm$: for the index $+$ the segment $[h_2,1]+m\tau+i0$
on the lower boundary of the smaller rectangle is identified with the lower bank $[h_2,1]+m\tau-i0$ of the slit in the larger rectangle. For the  index $-$ the segment $[h_2,1]+(m+1)\tau-i0$ on the upper boundary of the smaller rectangle is identified with the upper bank $[h_2,1]+(m+1)\tau+i0$ of the slit (see left panels of Fig. \ref{Sti3Graph}). 

Both decagons are topological discs and may be  conformally mapped to the upper half plane by functions $x_{3\pm}(u|\tau;n,m;h_1,h_2)$ uniquely defined up to the left action of  $PSL_2(\mathbb{R})$.
Our routine argument involving the reflection principle suggests that the function 
\be
\label{Sti3pm}
St_{3\pm}(\cdot|\tau;n,m;h_1,h_2):= x(\cdot|\tau)\circ x_{3\pm}^{-1}(\cdot|\tau;n,m;h_1,h_2). 
\ee
is a degree $n$ real rational function. Its qualitative graph is shown on the right panel of the Fig. \ref{Sti3Graph}. Note that $St_{3\pm}$ has $n-1$ full oscillations at the stopband $E^-$  and exactly $m+1$ and $n-m-1$ times on two  passbands $E_1^+$ and $E^+_2$ defined as follows: 
\be
\begin{array}{l}
E^-=x_{3\pm}(u|\tau; n, m; h_1,h_2),\qquad  u\in-1+[0,n-1]\tau;\\
E_1^+=x_{3\pm}(u|\tau; n, m; h_1,h_2),\qquad  u\in 1+[0,m+1]\tau;\\
E_2^+=x_{3\pm}(u|\tau; n, m; h_1,h_2),\qquad  u\in 1+[m,(n-1)]\tau,
\end{array}
\ee
 After the suitable rescaling functions $St_{3\pm}(x,|)$ acquire $2n+2$ alternation points and therefore become solutions of appropriate three band filter problems.
Both of them have topological class $\sigma=$ $(1,1,n)~mod~2$.

The auxiliary function $x_3^\pm(u|\dots)$ conformally mapping the decagon to the upper half plane is the inversion of the holomorphic abelian integral living on the genus three hyperelliptic curve with purely real branchpoints. Six  of the latter are the endpoints of the working bands $E$ and two more appear in the gap between the passbads $E_1^+$ and $E_2^+$. The latter two branchpoints bound  a 'fake' stopband with exactly one oscillation of the filtering function on it. Note that $St_{3\pm}$ is not monotonic away from the working bands and there are two real critical points in one gap between a passpand and the 'fake' stopband. Two branchpoints of the associated Riemann surface which are not the endpoints of $E$ are free parameters of the Ansatz along with the modular parameter $\tau$. Given $E$ one can write the set of three equations for those three parameters.

In the concluding sections of the paper we give a sketch of the prove that each solution of full degree $n$ has a representation of one of the mentioned 6 types as the composition of the elliptic sine and the holomorphic hyperelliptic  integral of genus $g=1,2$ or $3$.

\section{The ABC of Chebyshev Ansatz} 
Chebyshev Ansatz is a low-parametric representation \cite{B10} of a rational function $R$ which obeys the following restriction: the vast majority of its critical points are simple with the values in a four element set $\sf Q$, w.l.o.g. we assume that ${\sf Q}=\{\pm1, \pm 1/k(\tau)\}$, $k\neq0,\pm1,\infty$. The number of the exceptional critical points defines the so called extremality number of the function by the formula
\be
g(R):=1+\sum\limits_{x:~R(x)\not\in\sf Q} \ord~ dR(x)+
\sum\limits_{x:~R(x)\in\sf Q} [\frac12\ord~ dR(x)],
\label{gR}
\ee
here summation is taken over the whole Riemann sphere, $\ord~ dR(x)\ge0$ is the order of zero of the differential of $R$ at $x$ in local coordinates (say, a simple pole of $R$ gives $\ord =0$) and $[\cdot]$ is the integer part of a number. 

The Ansatz takes the appearance
\be
R(x)=sn(\int_{(e,0)}^{(x,w)}d\zeta+A(e)|\tau), 
\qquad A(e)= K(\tau)\left\{
\begin{array}{ll}
\pm 1, & R(e)=\pm1\\
\pm1+\tau, & R(e)=\pm 1/k(\tau).
\end{array}
\right.
\label{Ansatz}
\ee
where $d\zeta$ is a holomorphic abelian differential on the (unknown beforehand)
hyperelliptic curve $M=M(\sf E)$ with the affine model
$$
M({\sf E}) =\{(x,w)\in\mathbb{C}^2:\quad w^2=\prod_{e\in{\sf E}}(x-e)\}.
$$
The branchpoints of the curve $M(\sf E)$ are exactly those points $e$ on the sphere
where $R(x)$ takes values in $\sf Q$ with odd multiplicity (e.g. simple values).
The genus $g(M)$ of the curve $M$ equals to the extremality number of the rational function $R(x)$ that generated the curve \cite{B10}. The curve $M$ generated by a rational function $R$ in this framework is not arbitrary, for instance it covers the torus $M(\sf Q)$ with due branching. This implies significant restrictions on the moduli of the curve described e.g. in \cite{B10}.

The differential $d\zeta\in\Omega^1(M)$ in \eqref{Ansatz} is the pullback of the (suitably scaled) holomorphic differential on the torus $M(\sf Q)$ under the holomorphic  mapping $M({\sf E})\to M({\sf Q})$ induced by $R(x)$. The details of all constructions may be found in \cite{B10}. The divisor (i.e. positions and multiplicities of zeros) of $d\zeta$ may be easily reconstructed from the representation \eqref{Ansatz}: it projects to the points $x\in\mathbb{CP}^1$ which give the nontrivial input to the extremality number $g(R)$ of the function in formula \eqref{gR}. This observation gives us an alternative proof of the equality $g(M)=g(R)$.

The role of the phase shift $A(e)$ in Chebyshev Ansatz \eqref{Ansatz} is to switch between the solutions $\pm R(x)$, $\pm 1/(kR(x))$ which allows to interchange stop- and passbands as long as to modify their topological classes. For instance, to build solutions with $\sigma_1=0$ from functions $St_*$.

Despite the fact that very few elements in the Chebyshev representation \eqref{Ansatz}
are known, there is a big temptation to use for finding solutions $R(x)$  of filter approximation problems, since the extremality number $g(R)$ is usually small compared to its degree $n$. As we shall see, for the three band filter problem the value $g$ is at most three and we have to solve at most three transcendential equations for three (real) variables to determine all parameters of the Ansatz.

\section{Analysis of curves associated to solutions of 3-band problems}
The useful technique for the analysis of multuband uniform approximation problem is to consider the maximal extension of bands keeping the value of the approximation error intact. 
See e.g. \cite{B99} for multiband Chebyshev polynomials.

For any  solution  $R(x)$ of a filter problem with the workbands $E$ 
we consider the set of its extended bands $E_\sharp:= R^{-1}(R(E)) \supseteq E$ which generally is a collection of real analytic arcs in the complex plane.  The branchpoints of the associated curve $M$ are exactly the odd degree vertices of $E_\sharp$ considered as a graph.


\begin{thrm}
Given a full degree $n$ solution $R(x)$ of the 3-band filter problem $(E,\sigma)$, 
the extension $E_\sharp$ of the workband set is real and its intersection with the transition bands consists of a (possibly empty) segment in lacuna $T_{12}$ (see \eqref{AllBands} for the definition). 
\end{thrm}

\subsection{Discussion} 

There are exactly three possibilities for the location of the compliment $E_\sharp\setminus E$:

1. Lacuna between the passbands disappears, i.e. $E_\sharp=E\cup T_{12}$.
The extended workband set consists of two segments.  This is Zolotarev case, $g(M)=1$.

2. One endpoint of lacuna $T_{12}$ moves inside it, other being fixed. 
The extended workband set consists of three segments, hence $g(M)=2$.
One can show that there is exactly one exceptional critical point of the solution $R(x)$ which is located in the squeezed lacuna $T_{12}\setminus E_\sharp$. The abelian integral of the 
Chebyshev Ansatz maps the upper half plane to a rectangle with a slit. 
This is Stiefel case considered in Sect \ref{Sect:Sti}.

3. A segment of $E^\sharp$ appears amidst of the lacuna $T_{12}$. 
Now $E_\sharp$ consists of four segments and $g(M)=3$.
Further analysis shows that the solution $R(x)$  maps this segment 1-1
to one of two segments $R(E^\pm)$. 

$3^+$~
In case of sign $+$  a fake passband with the oscillation number $1$ appears in the lacuna $T_{12}$. Each of the segments in the compliment of the lacuna to the fake passband contains an
exceptional critical point of the solution $R(x)$. The abelian integral in \eqref{Ansatz}
maps the upper half plane to a twice slit rectangle we considered in Sect. \ref{Sect:Unival}

$3^-$~
In case of sign $-$  a fake stopband with the oscillation number $1$ appears in the lacuna between two passbands $E_1^+$ and $E_2^+$.
Deep analysis of combinatorics of the extended workband set shows that the solution $R(x)$ also has two exceptional critical points (which lift to the zeros of associated differential $d\zeta\in\Omega^1(M)$).
Those points may be either complex conjugate which give us the rectangular octagon of Sect. \ref{Sect:Octagon} with the internal branching. Or otherwise both of them lie in the same component of the  compliment of lacuna $T_{12}$ to the fake stopband which give us two rectangular decagons from the Sect. \ref{Sect:Decagon} as the image of $\mathbb{H}$ under the abelinan integral mapping.\\[3cm]

Topics that did not receive due attention in this  paper (e.g. proof of Theorem 1; the equations for the parameters of the Ansatz; etc.) will be considered in its extended version.


\begin{verbatim}
Institute for Numerical Math., Russian Academy of Sciences, 
Address: Russia 119991 Moscow ul. Gubkina, 8   
Email: ab.bogatyrev@gmail.com  
\end{verbatim}

\end{document}